\documentclass[12pt]{article}
\usepackage{amsmath,amsthm, amssymb}
\def\Omega{\varOmega}

\def\XXint#1#2#3{{\setbox0=\hbox{$#1{#2#3}{\int}$}
	\vcenter{\hbox{$#2#3$}}\kern-.5\wd0}}

\begin{document}

\title{\large 
\bf On the maximum principle for parabolic equations \\
with unbounded coefficients}
\author { {\it A.I.~Nazarov}\footnote {Partially supported
by St.Petersburg University grant 6.38.670.2013 and by the grant for support
of Leading Scientific Schools of Russia NSh-1771.2014.1.} ,\\ 
St.Petersburg Dept of Steklov Institute and St.Petersburg University, \\
{\small e-mail: \ al.il.nazarov@gmail.com}} 
\date{}\maketitle

\section{Introduction}

This text is based on the paper \cite{N87} and the note \cite{N88} published
in Russian in collected volumes by the Institute of Mathematics, Siberian Branch of USSR Academy of Sciences.
Later it turned out that the proofs in \cite{N87} can be essentially simplified. In particular, high-level
arguments from \cite{BL} and \cite{KP} can be avoided (see in this connection \cite{N01}). Also we fixed
some gaps in proofs of auxiliary assertions.\bigskip

We consider a priori maximum estimates for solution of initial-boundary value problem to parabolic equation
$$
{\cal L}u := \sigma(x,t) D_tu-a_{ij}(x,t)D_iD_ju+b_i(x,t)D_iu+c(x,t)u=f(x,t)	\eqno(1.1)
$$
in terms of the right-hand side in various spaces. Here and elsewhere we adopt the convention regarding 
summation with respect to repeated indices.

Such estimates for the Dirichlet problem to elliptic equations were established by A. D. Aleksandrov \cite{Al},
\cite{Al1}. N. V. Krylov \cite{Kr1}, \cite{Kr2} obtained these estimates for parabolic equations via
$\|f\|_{n+1,Q}$ provided all coefficients are bounded. N. N. Uraltseva and author \cite{NU} succeeded
to replace this assumption for $b_i$ by $b_i\in L_{n+1}(Q)$. Similar results were independently obtained by
Kai-sing Tso \cite{Ts} using a different method. Finally, N. V. Krylov \cite{Kr4} unified the estimates
of \cite{NU}, \cite{Ts}. Also he obtained the estimate via $\|f\|_{n+1,Q}$ provided 
$b_i\in L_n^xL_\infty^t(Q)$, and similar estimates via $\|f\|_{p+1,Q}$, $p\ge n$.

We establish the estimates of the same type in the space scale $L_p^xL_q^t$ (for $p\le q$) or $L_q^tL_p^x$ 
(for $p\ge q$) with arbitrary $p,q\le\infty$ subject to $\frac{n}{p}+\frac{1}{q}\le1$. Coefficients
$b_i$ are assumed to belong to a space of the same type, maybe with different $p$ and $q$. Moreover,
we can manage the ``composite'' coefficients
$$
b_i=\sum\limits_{k=1}^m b_i^{(k)}, \qquad b_i^{(k)}\in L_{p_k}L_{q_k}.	\eqno(1.2)
$$

The paper is organized as follows. Section 2 is devoted to basic estimates. In Section 3 we prove the pivotal 
lemma and then derive the required estimates in non-degenerate case. In Section 4 we generalize these
estimates for more wide class of operators. Also we prove the so-called Bony-type maximum principle. 
The estimate for operators with ``composite'' coefficients is proved in Section 5.\medskip

Let us recall some notation. $x=(x_1,\dots,x_n)$ is a vector 
in $\mathbb R^n$ with the Euclidean norm $|x|$; $(x;t)$ is a point 
in $\mathbb R^{n+1}$.

$B_R=\{x\,\big|\,|x|<R\}$ is a ball in $\mathbb R^n$. 

$C_R=B_R\times\mathbb R^1$; $C_{RT}=B_R\times]0,T[$.

$Q\subset C_{RT}$ is a domain in $\mathbb R^{n+1}$; $\Omega$ is the projection of $Q$ to $\mathbb R^n_x$;
$\overline Q$ is the closure of $Q$; $\chi_Q$ is the characteristic function of $Q$.

$|Q|$ and $|\Omega|$ stand for the Lebesgue measure of corresponding dimension.

$\partial Q$ is boundary of $Q$ while $\partial'Q$ is its parabolic boundary that is the set of 
$(x^0,t^0)\in\partial Q$ such that there exist $\delta>0$ and a function $x(t)\in {\cal C}(\mathbb R^1)$
satisfying $x(t^0)=x^0$ and $(x(t),t)\in Q$ for $t\in ]t^0,t^0+\delta]$. In particular, if $Q= C_{RT}$ then
$\partial Q=B_R\times\{0\}\cup \partial B_R\times [0,T[$.

By $\sup\limits_Q u$ we denote the essential supremum of a function $u$ on a set $Q$. If $u$ is continuous
then $Q_u=\{(x,t)\,\big|\,u>0\}$. 

The symbol $D_i$ denotes the operator of differentiation with respect to 
$x_i$; in particular, $Du=(D_1u, \dots, D_nu)$ is the gradient of 
$u$. $D_tu$ stands for the derivative of $u$ with respect to $t$.

We always assume that in (1.1) $\sigma\ge0$, $a_{ij}\lambda_i\lambda_j\ge0$ for $\lambda\in\mathbb R^n$, and
$c\ge0$. ${\bf Sp} (a)$ stands for the trace of the matrix $a=(a_{ij})$.

${\cal C}(\overline{Q})$ is the space of continuous functions with the norm $\|\cdot\|_Q$. 
${\cal C}_0(\overline{Q})$ is the subspace of ${\cal C}(\overline{Q})$ consisting of functions vanishing on
$\partial Q$. ${\cal C}^\infty(\overline{Q})$ is the set of smooth functions 
in $\overline{Q}$.

Let $p,q\ge1$ and let $w(x,t)>0$ a.e. in $Q$. We define $L_p^xL_q^t[w](Q)$ as the space of 
(equivalence classes of) functions $u$ such that the norm
$$
\|u\|=\bigg[\int\limits_\Omega dx\bigg[\int\limits_0^T|wu|^q\,dt\bigg]^{\frac pq}\,\bigg]^{\frac 1p}
$$
is finite ($u$ is assumed to be extended by zero on $C_{RT}\setminus Q$). If $p$ or $q$ is infinite then
corresponding integral should be replaced by $\sup$.
Analogously, $L_q^tL_p^x[w](Q)$ is the space with norm in which integrals are taken in reverse order.
If $w\equiv 1$ it is omitted.

By Minkowski's inequality, for $p<q$ the space $L_p^xL_q^t(Q)$ is continuously embedded into $L_q^tL_p^x(Q)$
(and $L_p^tL_q^x(Q)$ is continuously embedded into $L_q^xL_p^t(Q)$). For the sake of brevity we denote by
$\|\cdot\|_{p,q,(Q)}$ the norm in $L_p^xL_q^t(Q)$ if $p<q$, and the norm in $L_q^tL_p^x(Q)$ if $p>q$.
Thus, it always stands for the stronger norm, the first index corresponds to the spatial variables and
the second one -- to the time variable. For $p=q$ we evidently have $L_p^tL_p^x(Q)=L_p^xL_p^t(Q)=L_p(Q)$.

$W^{2,1}_{p,q}(Q)$ is the space with norm
$$
\|u\|_{W^{2,1}_{p,q}(Q)}=\|u\|_{p,q,(Q)}+\|D_tu\|_{p,q,(Q)}+\|Du\|_{p,q,(Q)}+\|D(Du)\|_{p,q,(Q)}.
$$

We set $f_+:=\max \{f,0\},\ \ f_-:=\max \{-f,0\}$ and denote by $p'$ the H\"older conjugate exponent for $p$. 
We use letters $M$, $N$ (with or without indices) to denote various constants. To indicate that, say, $N$ 
depends on some parameters, we list them in the parentheses: $N(\dots)$.

\section{Nondegenerate case. Basic estimates}

In Sections 2 and 3 we suppose that 
$$\delta\le\sigma, c\le\delta^{-1};\qquad |b|\le\delta^{-1};\qquad
\delta|\lambda|^2\le a_{ij}\lambda_i\lambda_j\le\delta^{-1}|\lambda|^2,\quad \lambda\in\mathbb R^n
$$
for some $\delta>0$.

{\bf Lemma 2.1}. {\it Let non-negative functions $A,B\in W^{2,1}_{\infty}(Q)\cap {\cal C}(\overline Q)$ 
satisfy ${\cal L}A\ge{\bf Sp}(a)$, ${\cal L}B\ge|b|$ a.e. in $Q$. Then for all functions 
$u\in W^{2,1}_{n+1}(Q)\cap {\cal C}(\overline Q)$ such that $u\big|_{\partial'Q}\le0$ the following estimate
holds:}
$$
u\le N_1(n)(B^2+A)^{\frac n{2(n+1)}}\cdot
\bigg\|\frac {({\cal L}u)_+}{(\sigma\det(a))^{\frac 1{n+1}}}\bigg\|_{n+1,(Q)}.	\eqno(2.1)
$$
\begin{proof} 
This statement is a particular case of \cite[Lemma 1.1]{Kr4}.
\end{proof}

{\bf Lemma 2.2}. {\it Under assumptions of Lemma 2.1, for all functions 
$u\in W^{2,1}_{n,\infty}(Q)\cap {\cal C}(\overline Q)$ such that $u\big|_{\partial'Q}\le0$ 
the following estimate holds:}
$$
u\le N_2(n)(B^2+A)^{\frac 12}\cdot
\bigg\|\frac {({\cal L}u)_+}{(\det(a))^{\frac 1n}}\bigg\|_{n,\infty,(Q)}.	\eqno(2.2)
$$
\begin{proof} 
We follow the scheme of proof of \cite[Lemma 3.3]{Kr4}. Let 
$$
f(x)=\chi_{\Omega}\cdot\sup\limits_{t}\frac {({\cal L}u)_+}{(\det(a))^{\frac 1n}}.
$$
We introduce a sequence $f_k\in{\cal C}_0^\infty(\mathbb R^n)$ such that $f_k\ge0$, 
$\|f_k-f\|_{n,(\mathbb R^n)}\to0$ as $k\to\infty$.

By \cite[Theorem III.2.3]{Kr3}, for arbitrary $\beta>0$ there exist $\psi_k\in W^2_{\infty}(\mathbb R^n)$
such that
$$
|D\psi_k(x)|\le \psi_k(x)\beta^{\frac 12};\qquad 
0\le\psi_k(x)\le N_2(n)\beta^{-\frac 12}\|f_k\|_{n,(\mathbb R^n)},	\eqno(2.3)
$$
and for any non-negative matrix $(\alpha_{ij})$
$$
-\alpha_{ij}D_iD_j\psi_k+\beta\psi_k{\bf Sp}(\alpha)-f_k(\det(\alpha))^{\frac 1n}\ge0.	\eqno(2.4)
$$

Now we consider functions
$$
\xi_k=u-\psi_k-\|\psi_k\|_Q\cdot(\beta A+\beta^{\frac 12} B).
$$
It is evident that $\xi_k\in W^{2,1}_{n,\infty}(Q)$ and $\xi_k\big|_{\partial'Q}\le0$.
We claim that $\xi_k\le N\|({\cal L}\xi_k)_+\|_{n,\infty,(Q)}$ with $N$ independent on $\xi_k$.

Indeed, let first $\xi_k\in W^{2,1}_{\infty}(Q)$. We introduce functions 
$$\varphi_k=\delta^{-1}\chi_{\Omega}\cdot\sup\limits_{t}({\cal L}\xi_k)_+;\qquad
\widetilde\varphi_k\in{\cal C}_0^\infty(B_{R+\delta^{-2}}); \quad \widetilde\varphi_k>\varphi_k;
\quad \|\widetilde\varphi_k\|_{n,(\mathbb R^n)}\le 2\|\varphi_k\|_{n,(\Omega)}.
$$

Example VIII.2.2 in \cite{Kr3} shows that there exists a solution $v_k\le0$ of the boundary value problem
for the Monge--Amp\`ere equation
$$
\det(D(Dv_k))=\frac 1{n^n}\widetilde\varphi_k^n \quad\mbox{in}\quad B_{R+\delta^{-2}};\qquad 
v_k\big|_{\partial B_{R+\delta^{-2}}}=0.
$$
Moreover, $|v_k|\le M\|\widetilde\varphi_k\|_{n,(\mathbb R^n)}$ with $M$ independent on $\xi_k$. 

Since $v_k$ is convex, $|Dv_k(x)|\le \delta^2|v_k(x)|$ in $\Omega$. This implies
$$
-{\cal L}v_k\ge {\bf Sp} (a\cdot D(Dv_k))-|b|\cdot |Dv_k|+c\cdot|v_k|\stackrel{*}{\ge}
n\cdot(\det(a\cdot D(Dv_k)))^{\frac 1n}\ge \delta\widetilde\varphi_k\ge({\cal L}\xi_k)_+
$$
in $Q$ (* is the arithmetic-geometric means inequality).

Note that $\xi_k+v_k\le 0$ on $\partial'Q$. By the maximum principle (see, e.g., \cite[Lemma III.3.6]{Kr3})
we obtain $\xi_k\le|v_k|\le2\delta^{-1}M\|({\cal L}\xi_k)_+\|_{n,\infty,(Q)}$, and the claim follows. For
$\xi_k\in W^{2,1}_{n,\infty}(Q)$ we arrive at this estimate by approximation.

Inequalities (2.3) and (2.4) give
$$
{\cal L}\xi_k\le ({\cal L}u)_+ +a_{ij}D_iD_j\psi_k-b_iD_i\psi_k
-\|\psi_k\|_Q\cdot(\beta{\bf Sp}(a)+\beta^{\frac 12}|b|)\le({\cal L}u)_+-f_k(\det(a))^{\frac 1n},
$$
and therefore
$$
\xi_k\le N\|(({\cal L}u)_+-f_k(\det(a))^{\frac 1n})_+\|_{n,\infty,(Q)}\to 0
\quad\mbox{as}\quad k\to\infty.
$$
By (2.3) we have
$$
u\le (\beta A+\beta^{\frac 12} B+1)\cdot N_2(n)\beta^{-\frac 12}\lim\limits_k\|f_k\|_{n,(\mathbb R^n)}.
$$
Finally, we minimize over $\beta$, and the Lemma follows.
\end{proof}

{\bf Remark 2.1}. The norms in the right-hand side of (2.1) and (2.2) can be taken over the set $Q_u$.
To prove it we can apply these estimates to $Q_u$ instead of $Q$.\medskip

{\bf Lemma 2.3}. {\it For all functions 
$u\in W^{2,1}_{\infty}(Q)\cap {\cal C}(\overline Q)$ such that $u\big|_{\partial'Q}\le0$ 
the following estimate holds:}
$$
u\le \bigg\|\frac {({\cal L}u)_+}{c}\bigg\|_{\infty,(Q_u)}.	\eqno(2.5)
$$
\begin{proof} 
For $w\equiv\big\|\frac {({\cal L}u)_+}{c}\big\|_{\infty,(Q_u)}$ we have ${\cal L}(u-w)\le0$ in $Q_u$
and $u\le w$ on $\partial' Q_u$. By the maximum principle we get (2.5).
\end{proof}

{\bf Lemma 2.4}. {\it For all functions 
$u\in W^{2,1}_{\infty,1}(Q)\cap {\cal C}(\overline Q)$ such that $u\big|_{\partial'Q}\le0$ 
the following estimate holds:}
$$
u\le \bigg\|\frac {({\cal L}u)_+}{\sigma}\bigg\|_{\infty,1,(Q_u)}.	\eqno(2.6)
$$
\begin{proof} 
Denote by $\Omega_u(\tau)$ the section of $Q_u$ by the plane $t=\tau$ and set
$$
w(t):=\int\limits_0^t\bigg\|\frac {({\cal L}u(\cdot,\tau))_+}{\sigma(\cdot,\tau)}
\bigg\|_{\infty,(\Omega_u(\tau))}\,d\tau.
$$ 
Then ${\cal L}(u-w)\le0$ in $Q_u$ and $u\le w$ on $\partial Q_u$. By the maximum principle we get 
$u\le \max w$, that gives (2.6).
\end{proof}

{\bf Remark 2.2}. All estimates in Lemmata 2.1--2.4 have the form $u\le M\|({\cal L}u)_+\|_{X(Q_u)}$. If 
$u\big|_{\partial'Q}=0$  then we can apply these estimates also to $-u$. This gives four estimates of the
form $|u|\le M\|{\cal L}u\|_{X(Q_u)}$.

\section{Nondegenerate case. Final estimates}

We recall that we denote by $\|\cdot\|_{p,q,(Q)}$ the norm in $L_q^tL_p^x(Q)$ if $p\ge q$ and the norm 
in $L_p^xL_q^t(Q)$ if $p\le q$. We also suppose that the assumptions from the beginning of Section 2 are 
fulfilled.\medskip

{\bf Pivotal Lemma}. {\it Let $\frac np + \frac 1q\le 1$, and let the functions $A$ and $B$ satisfy the assumptions
of Lemma 2.1. Then for all $u\in W^{2,1}_{p,q}(Q)\cap {\cal C}(\overline Q)$ such that 
$u\big|_{\partial'Q}\le0$ the following estimate holds:}
$$
u\le N(n)\|B^2+A\|_{Q_u}^{\frac n{2p}}\cdot
\bigg\|\frac {({\cal L}u)_+}{\sigma^{\frac 1q}(\det(a))^{\frac 1p}\,c^{1-\frac np-\frac 1q}}
\bigg\|_{p,q,(Q_u)}.	\eqno(3.1)
$$
\begin{proof} 
We prove (3.1) in several steps.\medskip

{\bf Step 1}. Suppose that $Q=C_{RT}$ and $a_{ij}$, $b_i$, $c$, 
$\sigma$ are smooth. Then for smooth functions $f$ the boundary value problem
$$
{\cal L}u=f\quad \mbox{in}\quad Q;\qquad u\big|_{\partial'Q}=0
$$
is uniquely solvable, see, e.g., \cite[Ch. 3]{F}. Denote this solution by $u={\cal L}^{-1}f$. Then 
${\cal L}^{-1}$ is evidently a linear operator from ${\cal C}^\infty(\overline{Q})$ to 
${\cal C}_0(\overline{Q})$.\medskip

{\bf 1a}. Let $\frac np + \frac 1q=1$ and $p<q<\infty$. Lemma 2.1 and Remark 2.2 show that ${\cal L}^{-1}$ 
can be extended to the operator from $L_{n+1}[(\sigma\,\det(a))^{-\frac 1{n+1}}](Q)$ to 
${\cal C}_0(\overline{Q})$, and 
$$
\|{\cal L}^{-1}\|\le M_1:=N_1(n)\|B^2+A\|_Q^{\frac n{2(n+1)}}. 
$$
Similarly, by Lemma 2.2 and Remark 2.2, ${\cal L}^{-1}$ can be extended to the operator from the closure 
of ${\cal C}^\infty(\overline{Q})$ in $L_n^xL_\infty^t[(\det(a))^{-\frac 1n}](Q)$ to 
${\cal C}_0(\overline{Q})$, and\footnote{Note that this closure coincides with the space
$L_n^x{\cal C}^t[(\det(a))^{-\frac 1n}](Q)$.} 
$$
\|{\cal L}^{-1}\|\le M_2:=N_2(n)\|B^2+A\|_Q^{\frac 12}.
$$

Consider adjoint operator ${\cal L}^{-1*}$ (with respect to duality $\langle u,v\rangle=\int\limits_Quv$).
It maps $L_1(Q)$ (as a closed subspace of $\big({\cal C}_0(\overline{Q})\big)'$) to 
$\big(L_{n+1}[(\sigma\,\det(a))^{-\frac 1{n+1}}](Q)\big)'
=L_{\frac {n+1}n}[(\sigma\,\det(a))^{\frac 1{n+1}}](Q)$. Furthermore, it maps also $L_1(Q)$ to 
$L_{\frac n{n-1}}^xL_1^t[(\det(a))^{\frac 1n}](Q)\subset\big(L_n^x{\cal C}^t[(\det(a))^{-\frac 1n}](Q)\big)'$,
since its image consists of functions. Its norms in these pairs do not exceed $M_1$ and $M_2$, respectively.

By the H\"older inequality,
$$
\Vert v\cdot \sigma^{\frac 1q}(\det(a))^{\frac 1p}\Vert_{L_{p'}^xL_{q'}^t(Q)}\le 
\Vert v\cdot (\sigma\,\det(a))^{\frac 1{n+1}}\Vert_{L_{\frac {n+1}n}(Q)}^{\theta}\cdot 
\Vert v\cdot (\det(a))^{\frac 1n}\Vert_{L_{\frac n{n-1}}^xL_1^t(Q)}^{1-\theta},
$$
where $\theta=\frac {n+1}q$. Therefore, ${\cal L}^{-1*}$ maps $L_1(Q)$ to 
$L_{p'}^xL_{q'}^t[\sigma^{\frac 1q}(\det(a))^{\frac 1p}](Q)$, and its norm does not exceed 
$M_1^\theta M_2^{1-\theta}$. This gives
$$
\|{\cal L}^{-1}f\|_Q\le N(n)\|B^2+A\|_Q^{\frac n{2p}}
\cdot\big\|\sigma^{-\frac 1q}(\det(a))^{-\frac 1p}f\big\|_{p,q,(Q)},	\eqno(3.2)
$$
where $N=\max\{N_1,N_2,1\}$.\medskip

{\bf 1b}. Let $\frac np + \frac 1q=1$ and $q<p<\infty$. By Lemma 2.4 and Remark 2.2, ${\cal L}^{-1}$ 
can be extended to the operator from $L_1^t{\cal C}^x[\sigma^{-1}](Q)$ to ${\cal C}_0(\overline{Q})$.
Turning to adjoint operator and interpolating between $L_{\frac {n+1}n}[(\sigma\,\det(a))^{\frac 1{n+1}}](Q)$
and $L_\infty^tL_1^x[\sigma](Q)$, we again arrive at (3.2).\medskip

{\bf 1c}. Let $\frac np + \frac 1q<1$, and $p,q<\infty$. We define $p_0=n+\frac pq$, $q_0=\frac {nq}p+1$,
such that $\frac n{p_0} + \frac 1{q_0}=1$ and $\frac {p_0}p=\frac {q_0}q$. By Lemma 2.3 and Remark 2.2, 
${\cal L}^{-1}$ can be extended to the operator from ${\cal C}[c^{-1}](\overline{Q})$ to 
${\cal C}_0(\overline{Q})$. From {\bf 1a} and {\bf 1b} one can see that it is continuous from the space
$$
\aligned
L_{p_0}^xL_{q_0}^t[\sigma^{-\frac 1{q_0}}(\det(a))^{-\frac 1{p_0}}](Q) \quad \mbox{for} \quad p\le q;\\
L_{q_0}^tL_{p_0}^x[\sigma^{-\frac 1{q_0}}(\det(a))^{-\frac 1{p_0}}](Q) \quad \mbox{for} \quad p\ge q
\endaligned
$$
to ${\cal C}_0(\overline{Q})$. Turning to adjoint operator and interpolating, we arrive at
$$
\|{\cal L}^{-1}f\|_Q\le N(n)\|B^2+A\|_Q^{\frac n{2p}}\cdot
\big\|\sigma^{-\frac 1q}(\det(a))^{-\frac 1p}\,c^{\frac np+\frac 1q-1}f\big\|_{p,q,(Q)}.	\eqno(3.3)
$$

{\bf Step 2}. Let $u$ be a smooth function, $u\big|_{\partial'Q}\le0$, and $p,q<\infty$.\medskip

{\bf 2a}. Suppose that $Q$ and coefficients of operator are as in Step 1. We define $u_1$ and $u_2$ as 
solutions of boundary value problems
\begin{eqnarray*}
{\cal L}u_1=({\cal L}u)_+ & \mbox{in} \quad Q; & u_1\big|_{\partial'Q}=0;\\
{\cal L}u_2=-({\cal L}u)_- & \mbox{in} \quad Q; & u_2\big|_{\partial'Q}=u\big|_{\partial'Q}.
\end{eqnarray*}
By the maximum principle $u_2\le0$. Applying (3.2) or (3.3) to $u_1$, we obtain (3.1) with $Q_u$ replaced
by $Q$.\medskip

{\bf 2b}. Suppose that $Q_u$ does not touch $\partial'C_{RT}$. We introduce a domain $\widetilde Q$ with
piecewise smooth boundary such that $Q_u\subset \widetilde Q\subset C_{RT}$. Then we consider a sequence
of Lipschitz functions $\zeta_k$ such that $\zeta_k=({\cal L}u)_+$ in $\widetilde Q$ and 
$\zeta_k\downarrow({\cal L}u)_+\cdot\chi_{\widetilde Q}$.

Denote by $u_k$ the solution of boundary value problem
$$
{\cal L}u_k=\zeta_k\quad \mbox{in} \quad C_{RT}; \qquad u_k\bigr|_{\partial' C_{RT}}= 0.
$$
Then evidently  $u_k\ge 0\ge u$ on $\partial' \widetilde Q$, and ${\cal L}u_k\ge{\cal L}u$ in $\widetilde Q$. 
By the maximum principle $u_k\ge u$ in $\widetilde Q$. We apply to $u_k$ in $C_{RT}$ the estimate obtained
in {\bf 2a} and pass to the limit as $k\to\infty$. It gives us (3.1) with $Q^u$ replaced by 
$\widetilde Q$.\medskip

{\bf 2c}. Since $p,q<\infty$, we can extend this estimate to arbitrary admissible coefficients and functions
$u$ by approximation.\medskip

{\bf 2d}. For arbitrary $Q_u$ we can consider functions $u_\varepsilon=u-\varepsilon$ and approximate 
$Q^u$ by domains $Q_{u_\varepsilon}\subset\widetilde Q_k\subset Q^u$ described in {\bf 2b}. Then we apply to
$u_\varepsilon$ in $\widetilde Q_k$ the estimate obtained in {\bf 2c}. Passage to the limit as $k\to\infty$
and then as $\varepsilon \to 0$ gives (3.1) in required form. The statement for $p,q<\infty$ is 
proved.\medskip

{\bf 3}. The cases $p=n$, $q=1$ and $p=q=\infty$ are considered in Lemmata 2.2, 2.4 and 2.3, 
respectively.\medskip

{\bf 3a}. Let $p=\infty$, $1<q<\infty$. Then we consider the estimate (3.1) for $\max\{q,nq'\}\le p<\infty$.
Since $\|\varphi\|_{p,q,(Q_u)}\le\|\varphi\|_{\infty,q,(Q_u)}\cdot|\Omega|^{\frac 1p}$, we obtain
$$
u\le N(n)\|B^2+A\|_{Q_u}^{\frac n{2p}}\cdot \delta^{-\frac {2n}p}\cdot|\Omega|^{\frac 1p}\cdot
\bigg\|\frac {({\cal L}u)_+}{\sigma^{\frac 1q}c^{1-\frac 1q}}\bigg\|_{\infty,q,(Q_u)}.
$$
Passage to the limit as $p\to\infty$ gives (3.1) for $p=\infty$.\medskip

{\bf 3b}. In a similar way, if $q=\infty$, $n<p<\infty$, then we consider the estimate (3.1) for 
large finite $q$ and pass to the limit using the embedding
$\|\varphi\|_{p,q,(Q_u)}\le\|\varphi\|_{p,\infty,(Q_u)}\cdot T^{\frac 1q}$.
\end{proof}

{\bf Remark 3.1}. For $p=q$ the estimate (3.1) was obtained by N.V. Krylov \cite{Kr4} in direct way. 
Interpolation method clarifies the nature of the weight $c^{\frac {n-p}{p+1}}$ in \cite{Kr4}. The result
of \cite[Lemma 3.3]{Kr4} for elliptic operators also can be obtained by interpolation between border spaces
$L_n[(\det(a))^{-\frac 1n}]$ and ${\cal C}[c^{-1}]$.\medskip

{\bf Corollary 3.1}. If there exists a function $B$ satisfying assumptions of Lemma 2.1 then the assertion of
Pivotal Lemma holds with $\|B^2+A\|_{Q_u}$ replaced by $(\|B\|_{Q_u}+R)^2$. This follows from
Lemma 1.2 in \cite{Kr4}.\medskip

{\bf Theorem 3.1}. {\it Let the assumptions in the beginning of Section 2 are satisfied. Suppose that 
$\frac n{p_0} + \frac 1{q_0}\le 1$ and $\frac n{p_1} + \frac 1{q_1}\le 1$. We put
$$
h=\sigma^{-\frac 1{q_1}}(\det(a))^{-\frac 1{p_1}}\,c^{\frac n{p_1}+\frac 1{q_1}-1}|b|.	\eqno(3.4)
$$
Then for all functions $u\in W^{2,1}_{p_0,q_0}(Q)\cap {\cal C}(\overline Q)$ such that 
$u\big|_{\partial'Q}\le0$ the following estimate holds:
$$
u\le M^{\frac n{p_0}}\cdot
\bigg\|\frac {({\cal L}u)_+}{\sigma^{\frac 1{q_0}}(\det(a))^{\frac 1{p_0}}\,c^{1-\frac n{p_0}-\frac 1{q_0}}}
\bigg\|_{p_0,q_0,(Q_u)},	\eqno(3.5)
$$
where $M$ depends only on $n$, $R$, $p_1$ and the norm $\|h\|_{p_1,q_1,(Q_u)}$.}

\begin{proof}
{\bf 1}. Let $p_0,q_0,p_1,q_1<\infty$. Then it is sufficient to obtain the estimate (3.5) for smooth 
coefficients and functions $u$ and then to pass to the limit. Moreover, we can assume that $Q_u$ does not 
touch $\partial'C_{RT}$.\medskip

As in the proof of Pivotal Lemma, Step {\bf 2b}, we approximate $Q_u$ by a domain $\widetilde Q$ with
piecewise smooth boundary such that $Q_u\subset \widetilde Q\subset C_{RT}$. Then we introduce a
sequence of operators ${\cal L}_k$ with smooth coefficients, satisfying the assumptions of Theorem,
such that ${\cal L}_k={\cal L}$ in $\widetilde Q$ and $|b^{(k)}|\downarrow|b|\cdot\chi_{\widetilde Q}$.

Denote by $B_k$ the solution of boundary value problem
$$
{\cal L}_kB_k=|b^{(k)}| \quad \mbox{in}\quad C_{RT}; \qquad B_k\bigr|_{\partial' C_{RT}}= 0.
$$
This function satisfies the assumptions of Lemma 2.1 for the operator ${\cal L}_k$. Therefore, we can apply
the estimate (3.1) with $p=p_1$, $q=q_1$, with regard to Corollary 3.1, to $u=\pm B_k$. This gives
$$
\|B_k\|_{C_{RT}}\le N(n)\cdot \big(\|B_k\|_{C_{RT}}+R\big)^{\frac n{p_1}}\|h_k\|_{p_1,q_1,({C_{RT}})}
$$
(here $h_k$ is defined by (3.4) with $b$ replaced by $b_k$).

Note that $q_1<\infty$ implies $p_1>n$. Therefore, if $\|B_k\|_{C_{RT}}>R$ then
$$
\|B_k\|_{C_{RT}}\le \Big[2^nN(n)^{p_1}\|h_k\|_{p_1,q_1,({C_{RT}})}^{p_1}\Big]^{\frac 1{p_1-n}}.	\eqno(3.6)
$$
We substitute this estimate to (3.1), pass to the limit as $k\to\infty$ and obtain the inequality (3.5)
with $M=N(n)\big(2R+[2^nN(n)^{p_1}\|h\|_{p_1,q_1,(\widetilde Q)}^{p_1}]^{\frac 1{p_1-n}}\big)$. Then we
finish the proof as in Step {\bf 2d} of the proof of Pivotal Lemma.\medskip

{\bf 2}. The estimate (3.5) for $p_0=\infty$ evidently follows from Lemma 2.3 ($q_0=\infty$), Lemma 2.4 
($q_0=1$) and Step {\bf 3a} in Pivotal Lemma ($1<q_0<\infty$).\medskip

{\bf 3}. Let $q_0=\infty$, $n<p_0<\infty$, and/or $q_1=\infty$, $n<p_1<\infty$. Then, as in Step {\bf 3b} in 
Pivotal Lemma, we can consider the estimate (3.5) for large finite $q_0$ ($q_1$), use the embedding theorem 
and pass to the limit as $q_0\to\infty$ ($q_1\to\infty$).\medskip

{\bf 4}. Let $1<q_1<\infty$, $p_1=\infty$, $p_0>n$. Using the estimate (3.5) for large finite $p_1$ and the 
embedding theorem, we arrive at
\begin{eqnarray*}
u&\le & N(n)\bigg(2R+\Big[2^nN(n)^{p_1}
\bigg\|\frac {\sigma^{-\frac 1{q_1}}c^{\frac 1{q_1}-1}|b|}{(\det(a))^{\frac 1{p_1}}\,c^{-\frac n{p_1}}}\bigg\|
_{\infty,q_1,(Q_u)}^{p_1}|\Omega|\Big]^{\frac 1{p_1-n}}\bigg)^{\frac n{p_0}}\\
&\times&\bigg\|\frac {({\cal L}u)_+}{\sigma^{\frac 1{q_0}}(\det(a))^{\frac 1{p_0}}\,
c^{1-\frac n{p_0}-\frac 1{q_0}}}\bigg\|_{p_0,q_0,(Q_u)}.
\end{eqnarray*}
The expression in large brackets does not exceed
$$
2R+2^{\frac n{p_1-n}}N(n)^{1+\frac n{p_1-n}}\delta^{-\frac {2n}{p_1-n}}
\big\|\sigma^{-\frac 1{q_1}}c^{\frac 1{q_1}-1}|b|\big\|
_{\infty,q_1,(Q_u)}^{1+\frac n{p_1-n}}|\Omega|^{\frac 1{p_1-n}}.
$$
We push $p_1\to\infty$ and obtain
$$u\le N(n)\bigg(2R+N(n)\big\|\sigma^{-\frac 1{q_1}}c^{\frac 1{q_1}-1}|b|\big\|_{\infty,q_1,(Q_u)}
\bigg)^{\frac n{p_0}}\cdot\bigg\|\frac {({\cal L}u)_+}{\sigma^{\frac 1{q_0}}(\det(a))^{\frac 1{p_0}}\,
c^{1-\frac n{p_0}-\frac 1{q_0}}}\bigg\|_{p_0,q_0,(Q_u)}.	\eqno(3.7)
$$
Then, as in part {\bf 3}, we derive the desired estimate for $p_1=q_1=\infty$, $p_0>n$.\medskip

{\bf 5}. Now let $q_1=1$, $p_0>n$. Since for $q>1$ and $\varphi\in L_{\infty}(Q_u)$ we have
$$
\|\varphi\|_{\infty,q,(Q_u)}\le
\|\varphi\|_{\infty,(Q_u)}^{\frac {q-1}q}\cdot\|\varphi\|_{\infty,1,(Q_u)}^{\frac 1q},
$$
the expression in brackets in (3.7) does not exceed
$$2R+N(n)\delta^{-4\frac {q_1-1}{q_1}}\big\|\sigma^{-1}|b|\big\|_{\infty,1,(Q_u)}^{\frac 1{q_1}}.
$$
We push $q_1\to1$ and obtain (3.7) for $q_1=1$.\medskip

In a similar way we consider the case $p_0=n$, $p_1>n$. For $u\in W^{2,1}_{\infty}(Q_u)$ we have from
part~{\bf 3}
$$
u\le M^{\frac n{p_0}}\cdot
\bigg\|\frac {({\cal L}u)_+}{(\det(a))^{\frac 1{p_0}}\,c^{1-\frac n{p_0}}}\bigg\|_{n,\infty,(Q_u)}^{\frac n{p_0}}\cdot
\bigg\|\frac {({\cal L}u)_+}{(\det(a))^{\frac 1{p_0}}\,c^{1-\frac n{p_0}}}\bigg\|_{\infty,(Q_u)}^{1-\frac n{p_0}}.
$$
Passage to the limit as $p_0\to n$ gives the desired estimate, and it remains to recall that
$W^{2,1}_{\infty}(Q_u)$ is dense in $W^{2,1}_{n,\infty}(Q_u)$.\medskip

{\bf 6}. The case $p_1=n$ is a special one since the inequality (3.6) fails. We construct a function $B$ from 
Pivotal Lemma in a different way, see \cite[Section 3]{Kr4}. We introduce a function 
$$f\in{\cal C}_0^\infty(B_{R+\varepsilon}); \quad f(x)>\sup\limits_t \,(h(x,t)\cdot\chi_{\Omega_u(t)});
\quad \|f\|_{n,(\mathbb R^n)}\le 2\|h\|_{n,\infty,(Q_u)}.
$$
Set $B:=-v$ where $v$ is the solution of boundary value problem
$$
\det(D(Dv))=\Big(\frac 2n\Big)^n f^n (1+|Dv|^2)^{\frac n2}\quad\mbox{in}\quad B_{R+\varepsilon};\qquad 
v\big|_{\partial B_{R+\varepsilon}}=0.
$$
Lemmata 3.1 and 3.2 in \cite{Kr4} and Remark 3.1 in \cite{Kr4} show that $B$ satisfies the assumptions of
Lemma 2.1, and
$$
\|B\|_{C_{RT}}\le N_3(n) (R+\varepsilon)\exp(N_4(n)\|f\|_{n,(\mathbb R^n)}^n).
$$
Finally, we can push $\varepsilon\to0$.
\end{proof}

{\bf Remark 3.2}. As it is pointed in Introduction, Theorem 3.1 and more general Theorem 4.1 were proved
by various methods for $p_0=q_0=n+1$, $p_1=q_1=\infty$ (see \cite{Kr2}); for $p_0=q_0=p_1=q_1=n+1$
(see \cite{NU}); for $p_0=q_0$, $p_1=q_1$ or $p_1=n$ (see \cite{Kr4}). See also \cite{Al}, \cite{Al1} 
for the case $p_0=p_1=n$.

\section{Generalization of Theorem 3.1}

In this Section we weaken requirements for coefficients of the operator ${\cal L}$ comparing to Sections 2 and 
3.\medskip

{\bf Theorem 4.1}. {\it Let $\frac n{p_0} + \frac 1{q_0}\le 1$ and $\frac n{p_1} + \frac 1{q_1}\le 1$. Suppose 
that the following assumption (depending on $p_0$ and $q_0$) is satisfied a.e. in $Q$:
\begin{eqnarray*}
 {\bf Sp}(a)>0 & \mbox{if} & p_0=n;\\
 \sigma >0     &  \mbox{if} & q_0=1;\\
 c>0           & \mbox{if} & p_0=q_0=\infty;\\
 c+\sigma>0    & \mbox{if} & p_0=\infty,\ 1<q_0<\infty;\qquad\qquad\qquad\qquad (4.1)\\
 {\bf Sp}(a)+c>0 & \mbox{if} & q_0=\infty,\ n<q_0<\infty;\\
 {\bf Sp}(a)+\sigma>0 & \mbox{if} & \frac n{p_0}+\frac 1{q_0}=1,\ p_0,q_0<\infty;\\
 \qquad\qquad\qquad\qquad\qquad\quad {\bf Sp}(a)+\sigma+c>0 & \mbox{if} & \mbox{otherwise}.
 \end{eqnarray*}
Let also $\|h\|_{p_1,q_1,(Q)}<\infty$, where the function $h$ is defined in (3.4). Then for all functions 
$u\in W^{2,1}_{p_0,q_0}(Q)\cap {\cal C}(\overline Q)$ such that $u\big|_{\partial'Q}\le0$, the estimate 
(3.5) holds. The quantity $M$ in (3.5) depends only on $n$, $R$, $p_1$ and the norm $\|h\|_{p_1,q_1,(Q_u)}$, 
and we set $0^0=1$, $\frac 00=0$, if such expression arises.}\medskip

\begin{proof}
{\bf 1}. Let  $p_0,q_0<\infty$ and $\frac n{p_0}+\frac 1{q_0}<1$. We set 
$$
{\cal L}_su:=\chi_{h\le s}\cdot{\cal L}u+\chi_{h>s}\cdot\big(D_tu-\Delta u+u\big).
$$
Let $a_{ijs}$, $b_{is}$, $c_s$, $\sigma_s$ be the coefficients of ${\cal L}_s$. Then ${\cal L}_s$ evidently 
satisfies assumptions of Theorem 4.1 with $h_s=h\cdot\chi_{h\le s}$, and
\begin{eqnarray*}
&&\bigg\|\frac {({\cal L}_su)_+}{\sigma_s^{\frac 1{q_0}}(\det(a_s))^{\frac 1{p_0}}
\,c_s^{1-\frac n{p_0}-\frac 1{q_0}}}\bigg\|_{p_0,q_0,(Q_u)}\\
\qquad\quad&\le& 
\bigg\|\frac {({\cal L}u)_+\cdot\chi_{h\le s}}
{\sigma^{\frac 1{q_0}}(\det(a))^{\frac 1{p_0}}\,c^{1-\frac n{p_0}-\frac 1{q_0}}}\bigg\|_{p_0,q_0,(Q_u)}+
\big\|(D_tu-\Delta u+u)_+\cdot\chi_{h> s}\big\|_{p_0,q_0,(Q_u)}.\qquad(4.2)
\end{eqnarray*}
Since $p_0,q_0<\infty$, the right-hand side of (4.2) tends to the norm in the right-hand side of (3.5) as 
$s\to\infty$. Thus, in this case it is sufficient to prove Theorem for $h$ bounded.\medskip

It  is evident that (3.5) does not change if we multiply all coefficients of ${\cal L}$ by the same function
positive almost everywhere. Thus, by (4.1) we can assume without loss of generality that 
${\bf Sp}(a)+\sigma+c=1$ a.e. in $Q$ and therefore all coefficients of ${\cal L}$ are bounded.

For $\varepsilon>0$ we set 
$$
{\cal L}_\varepsilon u:={\cal L}u+\varepsilon\cdot\big(D_tu-\Delta u+u\big).
$$
The operator ${\cal L}_\varepsilon$ satisfies all assumptions of Theorem 3.1, and $h_\varepsilon\le h$.
Therefore, the estimate (3.5) holds for ${\cal L}_\varepsilon$ instead of ${\cal L}$.

It remains to push $\varepsilon$ to zero and to note that ${\bf Sp}(a)+\sigma+c=1$ a.e. in $Q$ implies
$$
\bigg\|\frac {\varepsilon\cdot\big(D_tu-\Delta u+u\big)_+}
{(\sigma+\varepsilon)^{\frac 1{q_0}}(\det(a+\varepsilon I))^{\frac 1{p_0}}\,
(c+\varepsilon)^{1-\frac n{p_0}-\frac 1{q_0}}}\bigg\|_{p_0,q_0,(Q_u)}\le M(u)\varepsilon^\delta
\stackrel {\varepsilon\to0} \longrightarrow0
$$
(here $\delta=\min\{\frac 1{q_0},\frac 1{p_0},1-\frac n{p_0}-\frac 1{q_0}\}$ and $I$ stands for identity 
matrix).\medskip

{\bf 2}. In the case $\frac n{p_0}+\frac 1{q_0}=1$, $p_0,q_0<\infty$, repeating the first step of the 
part {\bf 1}, we reduce the proof to the case of bounded $h$ and ${\bf Sp}(a)+\sigma=1$ a.e. in $Q$.

For $s>0$, $\varepsilon>0$ we set $c_s=\min\{c,s\}$; 
$b_{is}=b_i\big(\frac {c_s}c\big)^{1-\frac n{p_1}-\frac 1{q_1}}$;
$$
{\cal L}_{s\varepsilon} u:=(\sigma+\varepsilon) D_tu-a_{ij}(x,t)D_iD_ju-\varepsilon\Delta u
+b_{is}D_iu+(c+\varepsilon)u.
$$
The operator ${\cal L}_{s\varepsilon}$ satisfies all assumptions of Theorem 3.1, and the estimate (3.5) 
holds for ${\cal L}_{s\varepsilon}$ instead of ${\cal L}$.

Since $p_0,q_0<\infty$, we can pass to the limit as $s\to\infty$. Then, similarly to part {\bf 1},
using ${\bf Sp}(a)+\sigma=1$ we push $\varepsilon$ to $0$.\medskip

{\bf 3}. For $p_0=\infty$, $1<q_0<\infty$ we can assume that $\sigma+c=1$ a.e. in $Q$. We apply the result of
part {\bf 1} to the operator ${\cal L}_\varepsilon$ for large finite $p$. By embedding 
$L_{q_0}^tL_\infty^x(Q)\to L_{q_0}^tL_p^x(Q)$ we have 
$$
u\le M^{\frac np}\varepsilon^{-\frac np}|\Omega|^{\frac 1p}\cdot\bigg\|\frac {({\cal L}_\varepsilon u)_+}
{(\sigma+\varepsilon)^{\frac 1{q_0}}\,(c+\varepsilon)^{1-\frac 1{q_0}}}\bigg\|_{\infty,q_0,(Q_u)}.
$$
We pass to the limit as $p\to\infty$. Then, similarly to part {\bf 1}, using $\sigma+c=1$ we push 
$\varepsilon$ to $0$.\medskip

In a similar way, for $q_0=\infty$, $n<p_0<\infty$ we can assume that ${\bf Sp}(a)+c=1$ a.e. in $Q$. 
We apply the result of part {\bf 1} to ${\cal L}_\varepsilon$ for large finite $q$ and obtain
$$
u\le M^{\frac n{p_0}}\varepsilon^{-\frac 1q}T^{\frac 1q}\cdot\bigg\|\frac {({\cal L}_\varepsilon u)_+}
{(\det(a+\varepsilon I))^{\frac 1{p_0}}\,(c+\varepsilon)^{1-\frac n{p_0}}}\bigg\|_{p_0,\infty,(Q_u)}.
$$
We pass to the limit as $q\to\infty$ and then as $\varepsilon\to0$.\medskip

In the same way, using these results we obtain the estimate for the case $p_0=q_0=\infty$.\medskip

{\bf 4}. Now let $p_0=n$. Then we can assume that ${\bf Sp}(a)=1$ a.e. in $Q$. For $p>n$, 
$u\in W^{2,1}_{\infty}(Q)\cap {\cal C}(\overline Q)$ we apply the result of part {\bf 1} to 
${\cal L}_\varepsilon$ and arrive at
$$
u\le M^{\frac np}\cdot\bigg\|\frac {({\cal L}_\varepsilon u)_+}
{(\det(a+\varepsilon I))^{\frac 1p}\,(c+\varepsilon)^{1-\frac np}}\bigg\|_{\infty,(Q_u)}^{1-\frac np}
\cdot\bigg\|\frac {({\cal L}_\varepsilon u)_+}
{(\det(a+\varepsilon I))^{\frac 1p}\,(c+\varepsilon)^{1-\frac np}}\bigg\|_{n,\infty,(Q_u)}^{\frac np}.
$$
Passing to the limit as $p\to n$ and then as $\varepsilon\to0$, we obtain the desired statement in this
case, since $W^{2,1}_{\infty}(Q)$ is dense in $W^{2,1}_{n,\infty}(Q)$.\medskip

The case $q_0=1$ is managed in a similar way.
\end{proof}

{\bf Remark 4.1}. For $p_0=\infty$ the constant in (3.5) does not depend on $h$. However, a simple example
shows that we cannot drop the restriction on $h$. Let
$$
Q=C_{1,1},\qquad {\cal L}u=D_tu-\Delta u+\frac {(n+1)x_i}{|x|^{\alpha}}D_iu+u.
$$
For $\alpha<2$ the operator ${\cal L}$ satisfies the assumptions of Theorem 4.1 since 
$\|h\|_{n,\infty,(Q)}<\infty$. However, if $\alpha=2$ then the function $U=2t-t^2-|x|^2-\frac 12$ satisfies
${\cal L}U<0$ while $U\big|_{\partial'Q}\le0$, $U(0,1)=\frac12$.\medskip

{\bf Remark 4.2}. The assumption $W^{2,1}_{p_0,q_0}(Q)$ in Theorem 4.1 can be replaced by 
$W^{2,1}_{p_0,q_0,loc}(Q)$. This fact can be proved as Lemma III.3.8 in \cite{Kr3}.\medskip

Now we weaken the assumption $c\ge0$. For the sake of brevity we formulate only the simplest 
generalization.\medskip

{\bf Theorem 4.2}. {\it Suppose that there is a constant $\varkappa>0$ such that 
$c_\varkappa:=c+\varkappa\sigma\ge0$ and the assumptions of Theorem 4.1 are satisfied with $c_\varkappa$ 
instead of $c$. Then for all functions $u\in W^{2,1}_{p_0,q_0}(Q)\cap {\cal C}(\overline Q)$ such that 
$u\big|_{\partial'Q}\le0$ the following estimate holds:}
$$
u\le \exp(\varkappa T)\cdot M^{\frac n{p_0}}\cdot
\bigg\|\frac {({\cal L}u)_+}{\sigma^{\frac 1{q_0}}(\det(a))^{\frac 1{p_0}}\,c^{1-\frac n{p_0}-\frac 1{q_0}}}
\bigg\|_{p_0,q_0,(Q_u)}.
$$

\begin{proof}
Consider the function $v=\exp(-\varkappa t)u$. We have 
${\cal L}_\varkappa v:={\cal L}v+\varkappa\sigma v=\exp(-\varkappa t){\cal L}u$. We apply Theorem 4.1 to
the operator ${\cal L}_\varkappa$ and to the function $v$. Then we take into account inequalities 
$\exp(-\varkappa t)\le1$ and $\exp(\varkappa t)\le\exp(\varkappa T)$, and the statement follows.
\end{proof}

Finally we prove the Bony-type maximum principle. In the case of bounded coefficients it was proved
in \cite{Bo} for elliptic operators and in \cite{Ts}, \cite{Kr4} for parabolic operators.\medskip

{\bf Theorem 4.3}. {\it Let the assumptions of Theorem 4.1 be satisfied with $p_0=p_1$, $q_0=q_1$. Suppose that a function 
$u\in W^{2,1}_{p_0,q_0,loc}(Q)$ attains its non-negative maximum in an interior point of $Q$. Then}
$$
\sup\limits_Q\frac {{\cal L}u}{{\bf Sp}(a)+\sigma+c}\ge0.
$$

\begin{proof}
Without loss of generality we can assume that ${\bf Sp}(a)+\sigma+c=1$. 

Let $\max\limits_Qu=u(x^0,t^0)\ge0$, $(x^0,t^0)\in Q$. Suppose that ${\cal L}u\le-\varepsilon<0$. Consider the cylinder 
$Q_\rho=\{(x,t)\,\big|\,t^0-\frac {\rho^2}2<t<t^0,\,|x|<\rho\}$ and introduce the function
$$
u_\delta(x,t)=u(x,t)-u(x^0,t^0)+\delta\Big(1-\frac {|x-x^0|^2-2(t-t^0)}{\rho^2}\Big).
$$
For sufficiently small $\rho$ we have $u\big|_{\partial Q_\rho}\le u(x^0,t^0)$ and therefore 
$u_\delta\big|_{\partial Q_\rho}\le0$. Thus, we can apply Theorem 4.1 with $p_0=p_1$, $q_0=q_1$. Since 
${\cal L}u_\delta\le-\varepsilon+\frac {\delta}{\rho^2}(2+2\rho|b|)$, this gives for 
$\delta<\frac {\varepsilon\rho^2}4$
$$
\delta=u_\delta(x^0,t^0)\le M(n,\rho,p_1,\|h\|_{p_1,q_1,(Q_\rho)})^{\frac n{p_1}}\cdot
\bigg\|\frac {\big(\frac {2\delta}{\rho}|b|-\frac {\varepsilon}2\big)_+}
{\sigma^{\frac 1{q_1}}(\det(a))^{\frac 1{p_1}}\,c^{1-\frac n{p_1}-\frac 1{q_1}}}\bigg\|_{p_1,q_1,(Q_\rho)}.
$$
Since 
$\sigma^{\frac 1{q_1}}(\det(a))^{\frac 1{p_1}}\,c^{1-\frac n{p_1}-\frac 1{q_1}}\le{\bf Sp}(a)+\sigma+c=1$,
we obtain
$$
\delta\le M^{\frac n{p_1}}\cdot\frac {2\delta}{\rho}\cdot
\big\|\big(h-\frac {\varepsilon\rho}{4\delta}\big)_+\big\|_{p_1,q_1,(Q_\rho)}=o(\delta)
\quad\mbox{as}\quad \delta\to0.
$$
This contradiction proves the statement.
\end{proof}

\section{The case of ``composite'' coefficients}

Consider the case where the coefficients $b_i$ are ``composite'', i.e. they can be written in the form (1.2),
and
$$
\|h_k\|_{p_k,q_k,(Q)}<\infty,\qquad
h_k=\sigma^{-\frac 1{q_k}}(\det(a))^{-\frac 1{p_k}}\,c^{\frac n{p_k}+\frac 1{q_k}-1}|b^{(k)}|	\eqno(5.1)
$$
for some $p_k,q_k$ such that $\frac n{p_k} + \frac 1{q_k}\le 1$. We again suppose that $c\ge0$ and other 
assumptions of Theorem 4.1.\medskip

{\bf Theorem 5.1}. {\it Under mentioned assumptions the estimate (3.5) holds with $M$ depending on
$n$, $R$, $p_1$ and norms $\|h_k\|_{p_k,q_k,(Q_u)}$, $k=1,\dots,m$}.

\begin{proof}
 We restrict ourselves to the case of smooth coefficients. Further arguments are similar to Section 4.\medskip
 
 Let $p_1=n$, and $p_k>n$ for $k\ge2$. Denote by $B_k$, $k\ge2$, solutions of boundary value problems
$$
{\cal L}_kB_k=|b^{(k)}| \quad \mbox{in}\quad C_{RT}; \qquad B_k\bigr|_{\partial' C_{RT}}= 0.
$$

As in part {\bf 6} of the proof ot Theorem 3.1, we introduce a function 
$$f\in{\cal C}_0^\infty(B_{R+1}); \quad f(x)>\sup\limits_t \,(h_1(x,t)\cdot\chi_{\Omega_u(t)});
\quad \|f\|_{n,(\mathbb R^n)}\le 2\|h_1\|_{n,\infty,(Q_u)}.
$$
Define $v$ as the solution of boundary value problem
$$
\det(D(Dv))=\Big(\frac 2n\Big)^n f^n (1+|Dv|^2)^{\frac n2}\quad\mbox{in}\quad B_{R+1};\qquad 
v\big|_{\partial B_{R+1}}=0
$$
and set
$$
B_1:=-v; \qquad \widetilde B:=\sum\limits_{k\ge2} B_k;\qquad B:=B_1+\widetilde B\cdot(1+\|DB_1\|_{C_{RT}}).
$$
Then
$$
{\cal L}B\ge(\det(a))^{\frac 1n}h_1(1+|DB_1|)-|b|\cdot|DB_1|+\sum\limits_{k\ge2} |b^{(k)}|
\cdot(1+\|DB_1\|_{C_{RT}})\ge|b|,
$$
and thus $B$ satisfies the assumptions of Pivotal Lemma.

We apply the estimate (3.1)  with $p=p_k$, $q=q_k$, with regard to Corollary 3.1, to $u=\pm B_k$. Summing
over $k\ge2$, we obtain
$$
\|\widetilde B\|_{C_{RT}}\le N(n)\sum\limits_{k\ge2} 
\big(\|B_1\|_{C_{RT}}+R+(1+\|DB_1\|_{C_{RT}})\|\widetilde B\|_{C_{RT}}\big)^{\frac n{p_k}}\cdot
\|h_k\|_{p_k,q_k,(C_{RT})},	\eqno(5.2)
$$
while Lemma 3.1 in \cite{Kr4} gives
$$
\|B_1\|_{C_{RT}},\ \|DB_1\|_{C_{RT}}\le N_3(n) (R+1)\exp(N_4(n)\|h_1\|_{n,\infty, (Q_u)}).
$$
Since $p_k>n$ for $k\ge2$, (5.2) easily implies the statement of Theorem.
\end{proof}

{\bf Remark 5.1}. The proof scheme of Theorem 5.1 is taken from \cite{Kr4}, where such proof was given in
some particular cases.\medskip

I am grateful to M.Z. Solomyak and N.N. Ural'tseva for important advices in 1987 when the original paper was
written. I also thank M.V. Safonov and S.V. Kislyakov who pointed me out possible simplification of
some proofs.


\begin{thebibliography}{AFT}

\bibitem[Al]{Al}
A.D.~Aleksandrov, {\em Uniqueness conditions and bounds for the solution of the Dirichlet problem},
Vestnik Leningrad. Univ. Ser. Mat. Meh. Astronom. {\bf 18} (1963), N3, 5--29 (Russian); English transl.: 
AMS Transl. (2) {\bf 68} (1968), 89--119.

\bibitem[Al1]{Al1}
A.D.~Aleksandrov, {\em Majorization of solutions of second-order linear equations}, 
Vestnik Leningrad. Univ. Ser. Mat. Meh. Astronom. {\bf 21} (1966), N1, 5--29 (Russian); English transl.: 
AMS Transl. (2) {\bf 68} (1968), 120--143.

\bibitem[BL]{BL}
 J.~Bergh, J.~L\"ofstr\"om, {\em Interpolation spaces. An introduction}, Springer, Berlin--Heidelberg, 1976.

\bibitem[Bo]{Bo} 
J.M.~Bony, {\em Principe du maximum dans les espaces de Sobolev}, C. R. Acad. sci. Paris {\bf 265} (1967),
N12, 333--336.
 
\bibitem[F]{F}
A. Friedman, {\em Partial differential equations of parabolic type}, Prentice-Hall (1964).
 
\bibitem[KP]{KP}
S.G.~Krein, Yu.I.~Petunin, {\em Scales of Banach spaces}, Uspekhi Mat. Nauk, {\bf 21} (1966), N2, 89--168
(Russian); English transl.: Russian Math. Surveys {\bf 21} (1966), N2, 85--159.
 
\bibitem[Kr1]{Kr1}
N.V.~Krylov, {\em Some estimates of the probability density of a stochastic integral},
Izv. AN SSSR, Ser. mat. {\bf 38} (1974), N1, 228--248 (Russian); English transl.:
Math. USSR -- Izvestiya, {\bf 8} (1974), N1, 233--254.

\bibitem[Kr2]{Kr2}
N.V.~Krylov, {\em Sequences of convex functions, and estimates of the maximum of the solution of 
a parabolic equation}, Sibirsk. Mat. Zh. {\bf 17} (1976), N2, 290--303 (Russian);
English transl.: Siberian Math. J. {\bf 17} (1976), N2, 226--236.

\bibitem[Kr3]{Kr3}
N.V.~Krylov, {\em Nonlinear elliptic and parabolic equations of the second order}, Moscow, Nauka, 1985 
(Russian); English transl.: Dordrecht, Reidel, 1987.

\bibitem[Kr4]{Kr4}
N.V.~Krylov, {\em On estimates for the maximum of solutions of a parabolic equation and estimates for
distribution of a semimartingal}, Mat. Sb. {\bf 130(172)} (1986), N2(6), 207--221 (Russian); 
English transl.: Math. USSR -- Sbornik, {\bf 58} (1987), N1, 207--221.

\bibitem[N87]{N87}
A.I.~Nazarov, {\em Interpolation of linear spaces and estimates for the maximum of 
a solution for parabolic equations}, Partial differential equations, Akad. 
Nauk SSSR Sibirsk. Otdel., Inst. Mat., No\-vo\-si\-birsk, 1987, 50--72 (Russian).

\bibitem[N88]{N88}
A.I.~Nazarov, {\em The maximum principle for parabolic equations with unbounded 
coefficients}, Some appli\-ca\-tions of functional analysis to problems of 
mathematical physics, Akad. Nauk SSSR Sibirsk. Otdel., Inst. Mat., 
No\-vo\-si\-birsk, 1988, 139--142 (Russian).

\bibitem[N01]{N01}
A.I.~Nazarov, {\em Estimates for the maximum of solutions of elliptic and 
parabolic equations in terms of weighted norms of the right-hand side}, 
Algebra \& Analysis {\bf 13} (2001), N2, 151--164 (Russian); English transl.: 
St.Petersburg Math. J. {\bf 13} (2002), N2, 269--279.

\bibitem[NU]{NU}
A.I.~Nazarov, N.N.~Uraltseva, {\em Convex-monotone hulls and an estimate of the 
maximum of the solution of a parabolic equation}, Zap. Nauchn. Sem. LOMI, 
{\bf 147} (1985), 95--109 (Russian); English transl.: J. Sov. Math. {\bf 37} 
(1987), N1, 851--859.

\bibitem[Ts]{Ts}
K.~Tso, {\em On an Aleksandrov-Bakel'man type maximum principle for 
second-order parabolic equations}, Comm. in PDE, {\bf 10} (1985), N5, 543--553.

\end{thebibliography}
\end{document}